\documentclass[a4paper,12pt]{article}

\usepackage{amsmath,amssymb,amsthm}
\usepackage{graphicx,color,hyperref}
\usepackage{rawfonts}

\newcounter{mathitem}
\newenvironment{mathitem}
  {\begin{list}{$(\roman{mathitem})$}{
   \setcounter{mathitem}{0}
   \usecounter{mathitem}
   \setlength{\topsep}{0pt plus 2pt minus 0pt}
   \setlength{\parskip}{0pt plus 2pt minus 0pt}
   \setlength{\partopsep}{0pt plus 2pt minus 0pt}
   \setlength{\parsep}{0pt plus 2pt minus 0pt}
   \setlength{\leftmargin}{35pt}
   \setlength{\itemsep}{0pt plus 2pt minus 0pt}}}
 {\end{list}}

\newtheorem{theorem}{Theorem}

\newtheorem{corollary}[theorem]{Corollary}

\title{The most general structure of graphs with hamiltonian or hamiltonian connected square}
\author{Jan Ekstein\thanks{Department of Mathematics and European Centre of Excellence NTIS - New Technologies for the Information Society, Faculty of Applied Sciences, University of West Bohemia, Pilsen, Technick\'a 8, 306 14 Plze\v n, Czech Republic, EU
\newline e-mail: \texttt{ekstein@kma.zcu.cz}.}\and
Herbert Fleischner \thanks{Institute of Logic and Computation, Algorithms and Complexity Group, Technical University of Vienna, Favoritenstrasse 9 - 11, 1040 Wien, Austria, EU  
\newline e-mail: \texttt{fleischner@ac.tuwien.ac.at}.}}
\date{\today}

\begin{document}
\maketitle

\begin{abstract}
On the basis of recent results on hamiltonicity, \cite{EkFle}, and hamiltonian connectedness, \cite{FleChia}, in the square of a 2-block, we determine the most general block-cutvertex structure a graph $G$ may have in order to guarantee that $G^2$ is hamiltonian, hamiltonian connected, respectively. Such an approach was already developed in \cite{FleHob} for hamiltonian total graphs. 
\medskip 
     
\medskip

\end{abstract}

\section{Introduction and Preliminary Discussion}
As for standard terminology and other terminology used in this paper, we refer to the book by Bondy and Murty, \cite{Bon}, and to the papers quoted in the references. Let $G$ be a connected graph. A \emph{2-block} is a 2-connected graph or a block of $G$ containing more than two vertices. The square of a graph $G$, denoted $G^2$, is the graph obtained from $G$ by joining any two nonadjacent vertices which have a common neighbor, by an edge.

It was shown in 1970 and published in 1974 that the square of every 2-block contains a       hamiltonian cycle, \cite{Fle2}. Key in proving this was the existence of EPS-graphs $S$ in connected bridgeless graphs $G$, where $S$ is the edge-disjoint union of a not necessarily connected eulerian subgraph $E$ and a linear forest $P$, and $S$ is connected and spans $G$, \cite{Fle1}. In subsequent papers \cite{Fle}, \cite{FleHob} the existence of various types of EPS-graphs was established. Their relevance was based on the fact that the total graph $T(G)$ of any connected graph $G$ other than $K_1$ is hamiltonian if and only if $G$ has an EPS-graph, \cite{FleHob}. This and the theory of EPS-graphs led to a description of the most general block-cutvertex graph $\mbox{bc}(G)$ of a graph $G$ may have such that $T(G)$ is hamiltonian and if $\mbox{bc}(G)$ does not have the corresponding structure, then exchanging certain 2-blocks in $G$ with some special 2-blocks yields a graph $G^*$ such that $\mbox{bc}(G)$ and $\mbox{bc}(G^*)$ are isomorphic but $T(G^*)$ is not hamiltonian, \cite{FleHob}. In dealing with hamiltonian cycles and hamiltonian paths by methods developed up to that point, it was shown in \cite{Fle} that in the square of graphs hamiltonicity and vertex-pancyclicity are equivalent concepts, and so are hamiltonian connectedness and panconnectedness. In this context Theorem~\ref{2-blockcycle} stated below was established as a tool needed to prove the equivalences just mentioned.

However, in the course of time much shorter proofs of Fleischner’s Theorem were developed \cite{Geo}, \cite{Riha}; the same applies to Theorem \ref{2-blockcycle} below, \cite{MutRau}. More recently, an algorithm yielding a hamiltonian cycle in the square of a 2-block in linear time, was developed, \cite{AlsGeo}. The methods developed in these much shorter proofs (including the algorithm just mentioned) do not seem to yield short proofs of Theorems \ref{H_4} and \ref{F_4} below, \cite{EkFle}, \cite{FleChia}. These latter theorems are, on the other hand, instrumental in proving the central results of this paper, i.e., Theorems \ref{hamiltonian} and \ref{hamconnected}, and related algorithms. 

Let $\mbox{bc}(G)$ denote the \emph{block-cutvertex graph} of $G$. Blocks corresponding to leaves of $\mbox{bc}(G)$ are called \emph{endblocks}, otherwise \emph{innerblocks}. Note that a block in a graph $G$ is either a 2-block or a bridge of $G$. For each cutvertex $i$ of $G$, let $k_{i}$ be the number of 2-blocks of $G$ which include vertex $i$ and let $\mbox{bn}(i)$ be the number of nontrivial bridges of $G$ which are incident with vertex $i$. In what follows a bridge is called nontrivial if it is not incident to a leaf.

Let $H$ be a subgraph of the graph $G$. We define $G-H:=G-E(H)-\{v\in V(H):d_H(v)=d_G(v)\}$. 
  
In Theorem \ref{hamiltonian}, we introduce an array $m_{i}(B)$ of numbers with an entry for each pair consisting of a cutvertex $i$ and a 2-block $B$ of $G$. We may think of this number $m_{i}(B)$ as the number of edges of $B$ incident with $i$ which are possibly contained in a hamiltonian cycle in $G^{2}$. 

Statement of Theorem \ref{hamiltonian} describes the most general block-cutvertex structure a graph $G$ may have in order to guarantee that $G^2$ is hamiltonian using parameters $m_{i}(B)$ as in \cite{FleHob}.

\begin{theorem}
 \label{hamiltonian}
 Let $G$ be a connected graph with at least three vertices. Let the 2-blocks of $G$ be labelled $B_{1},B_{2},...,B_{n}$. Let the cutvertices of $G$ be labelled $1,2,...,s$. Suppose there is a labelling $m_{i}(B_{t})$ for each $i\in\{1,2,...,s\}$ and each $t\in\{1,2,...,n\}$ such that the following conditions are fulfilled.
 \begin{mathitem}
  \item[1)] $0\leq m_{i}(B_t)\leq2$ for all $i$ and all 2-blocks $B_t$;
  \item[2)] for 2-block $B_t$ $m_{i}(B_t)=0$ if and only if cutvertex $i$ is not in $V(B_t)$;
  \item[3)] for 2-block $B_t$, $m_{i}(B_t)\geq\mbox{bn}(i)$, if cutvertex $i\in V(B_t)$;
  \item[4)] $\mbox{bn}(i)\leq2$ for all $i\in\{1,2,...,s\}$; 
  \item[5)] $\sum_{i=1}^{s}m_{i}(B_t)\leq4$ for each 2-block $B_t$ of $G$ and, if $m_i(B_t)=2$ 
            for some $i$, then $\sum_{i=1}^sm_{i}(B_t)\leq3$; and
  \item[6)] $\sum_{t=1}^{n}m_i(B_t)\geq2k_{i}+\mbox{bn}(i)-2$ for each $i\in\{1,2,...,s\}$.
 \end{mathitem}
 Then $G^{2}$ is hamiltonian. 

Moreover, if the labelling $m_{i}(B_{t})$  satisfying conditions 1), 2) and 3) is given and at least one of conditions  4), 5), 6) is violated by some $G$, then there exists a class of graphs $G'$ with non-hamiltonian square but $\mbox{bc}(G')$ and $\mbox{bc}(G)$ are isomorphic.
\end{theorem}

Also, we obtain a similar result for hamiltonian connectedness (Theorem~\ref{hamconnected}). Quite surprisingly, its formulation is much simpler than that of Theorem \ref{hamiltonian}.

\begin{theorem}
 \label{hamconnected}
 Let $G$ be a connected graph such that the following conditions are fulfilled:
 \begin{itemize}
  \item[1)] there is no nontrivial bridge of $G$;
  \item[2)] every block contains at most 2 cutvertices.
 \end{itemize}
 Then $G^2$ is hamiltonian connected.

Moreover,
\begin{itemize}
  \item[$\cdot$] if a graph $G$ contains a nontrivial bridge, then $G^2$ is not hamiltonian connected;  
  \item[$\cdot$] if $G$ contains a block containing more than 2 cutvertices, then there is a graph $G'$ such that $\mbox{bc}(G)$ and $\mbox{bc}(G')$ are isomorphic but $(G')^2$ is not hamiltonian connected. 
\end{itemize}
\end{theorem}

A fundamental result regarding hamiltonicity in the square of a 2-block is the following theorem.

\begin{theorem}\emph{\textbf{\cite{Fle}}}
\label{2-blockcycle}
Suppose $v$ and $w$ are two arbitrarily chosen vertices of a $2$-block $G$.
Then $G^2$ contains a hamiltonian cycle $C$ such that the edges of $C$ incident to $v$ are in $G$ and at least one of the edges of $C$ incident to $w$ is in $G$. Furthermore, if $v$ and $w$ are adjacent in $G$, then these are three different edges.  
\end{theorem}

The hamiltonian theme in the square of a 2-block has been recently revisited (\cite{EkChiaFle}, \cite{EkFle}, \cite{FleChia}), yielding the following results which are essential for this paper. 

A graph $G$ is said to have the \emph{$\mathcal{H}_{k}$ property} if for any given vertices $x_1,...,x_k$ there is a hamiltonian cycle in $G^2$ containing distinct edges $x_1y_1,...,x_ky_k$ of~$G$.
 
\begin{theorem}\emph{\textbf{\cite{EkFle}}}
 \label{H_4}
 Given a 2-block $G$ on at least 4 vertices, then $G$ has the $\mathcal{H}_{4}$ property, and there are 2-blocks of arbitrary order greater than 4 without the $\mathcal{H}_{5}$ property.
\end{theorem}

By a \emph{$uv$-path} we mean a path from $u$ to $v$ in $G$. If a $uv$-path is hamiltonian, we call it a \emph{$uv$-hamiltonian path}. Let $A=\{x_{1},x_{2},...,x_{k}\}$ be a set of $k\geq 3$ distinct vertices in $G$. An $x_{1}x_{2}$-hamiltonian path in $G^{2}$ which contains $k-2$ distinct edges $x_{i}y_{i}\in E(G), i=3,...,k$, is said to be $\mathcal{F}_{k}$. A graph $G$ is said to have the \emph{$\mathcal{F}_{k}$ property} if, for any set $A=\{x_{1},x_{2},...,x_{k}\}\subseteq V(G)$, there is an $\mathcal{F}_{k}$ $x_{1}x_{2}$-hamiltonian path in $G^{2}$.

\begin{theorem}\emph{\textbf{\cite{FleChia}}}
 \label{F_4}
  Every 2-block on at least 4 vertices has the $\mathcal{F}_{4}$ property. 
\end{theorem}

A graph $G$ is said to have the \emph{strong $\mathcal{F}_{3}$ property} if, for any set of 3 vertices $\{x_{1},x_{2},x_{3}\}$ in $G$, there is an $x_{1}x_{2}$-hamiltonian path in $G^{2}$ containing distinct edges $x_{3}z_{3},x_{i}z_{i}\in E(G)$ for a given $i\in\{1,2\}$. Such an $x_{1}x_{2}$-hamiltonian path in $G^{2}$ is called a strong          $\mathcal{F}_{3}$ $x_{1}x_{2}$-hamiltonian path.

\begin{theorem}\emph{\textbf{\cite{FleChia}}}
 \label{strongF_3}
  Every 2-block has the strong $\mathcal{F}_{3}$ property.
\end{theorem}

\begin{theorem}\emph{\textbf{\cite{FleChia}}}
	\label{strongF_3ends}
Let $G$ be a $2$-connected graph and let $x, y$ be two vertices in $G$. Then $G^2$ has an $xy$-hamiltonian path $P(x, y)$ such that

(i) $xz \in E(G) \cap E(P(x,y))$  for some $z \in V(G)$, 
and

(ii) either $yw \in  E(G) \cap E(P(x,y))$ for some $w\in V(G)$,  or  else $P(x, y)$ contains an edge $uv$ for some vertices $u, v \in N(y)$.
\end{theorem} 

\section{Proofs and algorithms}

PROOF OF THEOREM \ref{hamiltonian}

\begin{proof}
 Set $P_0=G-\cup_{t=1}^{n}B_t$. Then every component of $P_0$ is a tree. Since by 4) $\mbox{bn}(i)\leq 2$ every component of $P_0$ is even a caterpillar. 
 
 For every caterpillar $T$ of $P_0$ except $T=K_2$ we have the following observation which can be proved easily. 

\medskip

\noindent
\emph{Observation: Let $T$ be a caterpillar with at least three vertices and $P=x_1x_2...x_m$ be some longest path in $T$. Then $T^2$ contains a hamiltonian cycle containing edges $x_1x_2, x_{m-1}x_m$  and different edges $u_jv_j$, where $u_j,v_j\in N_G(x_j)$ for $j=2,3,...,m-1$.}

See Figure \ref{Catter} for illustration in which for $x_3$ we have $u_3=x_2$ and $v_3=x_4$.

\begin{figure}[ht]
\begin{center}
\includegraphics[width=9cm]{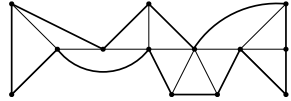}
\end{center}\caption{Hamiltonian cycle in a caterpillar for $m=7$ (bold edges)}\label{Catter}
\end{figure} 

\medskip

 Every 2-block $B_t$ contains a hamitonian cycle in $(B_t)^2$ which is one of two types depending on labellings $m_{i}(B_{t})$:
 
Let $m_i(B_t)\neq 2$ for every $i=1,2,...,s$. If $B_t\cong C_3$, then we set $C_t=B_t$. Otherwise for at most 4 cutvertices $a,b,c,d$ it holds that $m_j(B_t)=1$ for $j=a,b,c,d$ by condition 5). By Theorem \ref{H_4}, $(B_t)^2$ has a hamiltonian cycle $C_t$ containing 4 different edges $aa',bb',cc',dd'$ of~$B_t$. 
 
If $m_i(B_t)=2$ for some $i\in\{1,2,...,s\}$, then at most one cutvertex $a$ has $m_a(B_t)=1$ by condition 5). By Theorem \ref{2-blockcycle}, $(B_t)^2$ has a hamiltonian cycle $C_t$ containing 3 different edges $ii',ii'',aa'$ of $B_t$.   

The union of hamiltonian cycles $C_t$ in $(B_t)^2$, for $t=1,2,...,n$, hamiltonian cycles in the square of each catepillar (nontrivial component of $P_0$) and trivial components of $P_0$  is a connected spanning subgraph $S$ of $G^2$. 

We construct a hamiltonian cycle $C$ in $G^2$ from $S$ repeating step by step the following procedure for every cutvertex $i$ of $G$ with $m_{i}(B)\geq 1$ for some 2-block $B$. 

If $i$ does not exist, then $n=0$ and $G=P_0$ is a caterpillar. Hence $S$ is a hamiltonian cycle in $G^2$. Otherwise we join all hamiltonian cycles from $S$ containing $i$ together with trivial components of $P_0$ containing $i$ to one cycle in the following way.

\medskip

\noindent
First assume that $\mbox{bn}(i)=0$.

By condition 6) we have  $\sum_{t=1}^{n}m_i(B_t)\geq 2k_i-2$. Without loss of ge\-nerality for $k_i>1$ we may assume that $m_i(B_1)\geq 1$, $m_i(B_2)\geq 1$ and $m_i(B_3)=m_i(B_4)=\cdots=m_i(B_{k_i})=2$, where $m_i(B_t)$  corresponds to the number of edges of $B_t$ incident to $i$ in $C_t$. If $k_i=1$, then by condition 2) we have $m_i(B_1)\geq 1$.

We find a cycle $C^i$ on $\cup_{r=1}^{k_i}V(C_r)\cup L$, where $L$ is the set of all leaves incident to $i$, by appropriately replacing edges of $C_r\cap B_r$, $r=1,2,...,k_i$, incident to $i$ (guaranteed by definition of $m_i(B_t)$) with edges of $G^2$ joining vertices in different $C_r$ adjacent to $i$ and leaves adjacent to $i$. Note that we preserve properties given by $m_j(B_t)$ for all $j\neq i$.

\noindent
Now assume that $\mbox{bn}(i)=1$.

By condition 6) we have  $\sum_{t=1}^{n}m_i(B_t)\geq 2k_i+1-2=2k_i-1$. Without loss of generality we may assume that $m_i(B_1)\geq 1$ and $m_i(B_2)=m_i(B_3)=\cdots=m_i(B_{k_i})=2$, where $m_i(B_t)$  corresponds to the number of edges of $B_t$ incident to $i$ in $C_t$.
Let $T$ be the component of $P_0$ containing $i$.

If $T=K_2=ii'$, where $i'$ is also a cutvertex of $G$ with $m_{i'}(B)\geq 1$ ($T$ is a trivial component of $P_0$), then we find a cycle $C^i$ on $\cup_{r=1}^{k_i}V(C_r)\cup V(T)$ containing the edge $ii'$ by appropriately replacing edges of $C_r\cap B_r$, $r=1,2,...,k_i$, incident to $i$ (guaranteed by definition of $m_i(B_t)$) with edges of $G^2$ joining $i'$ and vertices in different $C_r$ adjacent to $i$. Also here we preserve properties given by $m_j(B_t)$ for all $j\neq i$.

If $T$ is a nontrivial component of $P_0$, then $T^2$ contains a hamiltonian cycle $C_T$ containing end-edges of any fixed longest path $P$ in $T$ (we choose end-edges containing cutvertices of $G$ with $m_{i}(B_t)\geq 1$) - see Observation above. Again we find a cycle $C^i$ on $\cup_{r=1}^{k_i}V(C_r)\cup V(C_T)$ by appropriately replacing edges of $C_r\cap B_r$, $r=1,2,...,k_i$, incident to $i$ (guaranteed by definition of $m_i(B_t)$) and the end-edge $ii^*$ of $P$ with edges of $G^2$ joining $i^*$ and vertices in different $C_r$ adjacent to $i$. Again we preserve properties given by $m_j(B_t)$ for all $j\neq i$ and by $C_T$.

\noindent
Finally assume that $\mbox{bn}(i)=2$.

By condition 6) we have  $\sum_{t=1}^{n}m_i(B_t)\geq 2k_i+2-2=2k_i$. It follows necessarily that  $m_i(B_1)=m_i(B_2)=\cdots=m_i(B_{k_i})=2$, where $m_i(B_t)$  corresponds to the number of edges of $B_t$ incident to $i$ in $C_t$.

Let $T$ be the nontrivial component of $P_0$ containing $i$. Note that $i$ is not an endvertex of $T$ because of $\mbox{bn}(i)=2$. Then $T^2$ contains a hamiltonian cycle $C_T$ containing end-edges of any fixed longest path in $T$ (we choose end-edges containing cutvertices of $G$ with $m_{i}(B_t)\geq 1$) and an edge $u_iv_i$ of $G^2$ where $u_i,v_i\in N_G(i)$  (see Observation above). We find a cycle $C^i$ on $\cup_{r=1}^{k_i}V(C_r)\cup V(C_T)$ by appropriately replacing edges of $C_r\cap B_r$, $r=1,2,...,k_i$, incident to $i$ (guaranteed by definition of $m_i(B_t)$) and the edge $u_iv_i$ of $P$ with edges of $G^2$ joining $u_i, v_i$ and vertices in different $C_r$ adjacent to $i$ if $k_i>1$. If, however, $k_i=1$, then $u_i$ and $v_i$ are joined to the neighbors of $C_r\cap B_r$ in $N_G(i)$. Also here we preserve properties given by $m_j(B_t)$ for all $j\neq i$ and by $C_T$.

Now we choose next cutvertex $i$ with $m_{i}(B)\geq 1$ for some 2-block $B$ successively and we use all cycles  formed in the previous steps instead of previously formed cycles. Note that we preserve all properties given by $m_j(B)$ for all $j\neq i$ in every case. We stop with the hamiltonian cycle in $G^2$ as required.

\vskip 8mm

Now assume that there is no labelling satisfying conditions 1) - 6), that is, the labelling $m_i(B_t)$ satisfying conditions 1), 2) and 3) is given and at least one of conditions  4), 5), 6) is violated. We show that there exists a class of graphs $G'$ with non-hamiltonian square but $\mbox{bc}(G')$ and $\mbox{bc}(G)$ are isomorphic.

\bigskip

\begin{figure}[ht]
\begin{center}
\includegraphics[width=13cm]{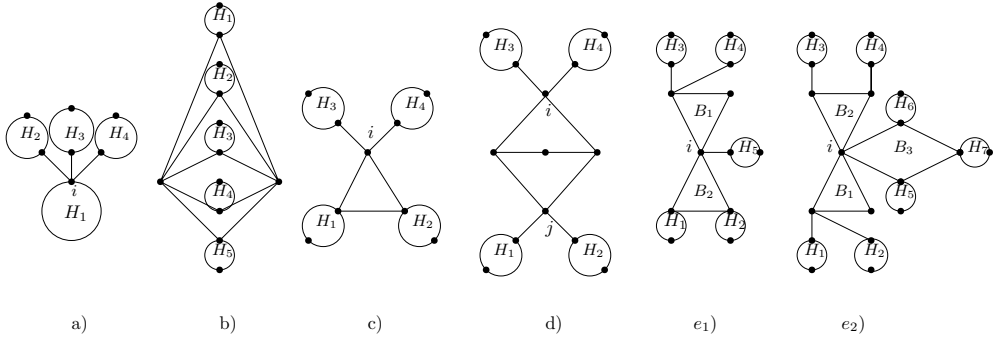}
\end{center}\caption{Graphs without hamiltonian square}\label{Counter}
\end{figure} 

\noindent
\emph{Condition 4) does not hold.} 

Hence $\mbox{bn}(i)\geq 3$ for at least one $i\in\{1,2,...,s\}$. Clearly this is a class of graphs $G'$ such that the square of every such graph $G'$ does not contain a hamiltonian cycle (if we try to construct a hamiltonian cycle in the square, then the degree of the cutvertex $i$ is at least 3, a contradiction), e.g. see the graph in Figure \ref{Counter} a), where $H_1$ is an arbitrary connected graphs, $H_2$, $H_3$, $H_4$ are arbitrary connected graphs with at least one edge each and $\mbox{bn}(i)=3$. Note that conditions 5) and 6) may hold.

\vskip 0.9cm

\noindent
\emph{Condition 5) does not hold.}

Hence $\sum_{i=1}^{s}m_{i}(B)\geq 5$ for some 2-block $B$ and $m_i(B)<2$ for all $i$ or $\sum_{j=1}^{s}m_{j}(B)\geq 4$ for some 2-block $B$ and $m_i(B)=2$ for some $i\in\{1,2,...,s\}$.

First suppose that $k=\sum_{i=1}^{s}m_{i}(B)\geq 5$ for some 2-block $B$ of $G$ and $m_i(B)<2$ for all $i$. Clearly $B$ has exactly $k$ cutvertices by condition 2). Then we exhange $B$ with $K_{2,k}$ where $k$ 2-valent vertices are cutvertices of $G$ and all other blocks with arbitrary blocks to get a class of graphs $G'$ such that $\mbox{bc}(G')$ and $\mbox{bc}(G)$ are isomorphic. The square of every such graph $G'$ does not contain a hamiltonian cycle (if we try to construct a hamiltonian cycle in the square, then the degree of at least one of the two $k$-valent vertices of $K_{2,k}$ is at least 3, a contradiction), e.g. see the graph in Figure \ref{Counter} b), where $k=5$ and $H_1,...,H_5$ are arbitrary connected graphs with at least one edge each. Note that conditions 4), 6) and the second part of condition 5) may hold.

Now suppose that $\sum_{j=1}^{s}m_{j}(B)\geq 4$ for some 2-block $B$ and $m_i(B)=2$ for some $i$. If $B$ contains at least 5 cutvertices of $G$, then we continue similarly as above. If $B$ contains $k$ cutvertices of $G$ where $2\leq k \leq 4$, then without loss of generality we may assume that we tried to set the labelling $m_i(B_t)$ satisfying firstly condition 5) and subsequently condition 6). Hence $\mbox{bn}(i)\geq 2$ and $\mbox{bn}(j)\geq 2$ where $j$ is the second cutvertex of $G$ in $B$ if $k=2$, otherwise we find a labelling $m_i(B_t)$ satisfying condition 5), a contradiction (see Algorithm~1 cases e) and f) below).

For $k=3,4$ we exhange $B$ with a cycle $C_k$ to get a class of graphs $G'$ such that $\mbox{bc}(G')$ and $\mbox{bc}(G)$ are isomorphic. The square of every such graph $G'$ does not contain a hamiltonian cycle (if we try to construct a hamiltonian cycle in the square, then the degree of the cutvertex $i$ is at least 3, a contradiction), e.g. see the graph in Figure \ref{Counter} c), where $k=3$ and $H_1,...,H_4$ are arbitrary connected graphs with at least one edge each. Note that conditions 4), 6) and the first part of condition 5) may hold.

For $k=2$, we exchange $B$ with $K_{2,3}$, where two of the three 2-valent vertices are $i$ and $j$, to get a class of graphs $G'$ such that $\mbox{bc}(G')$ and $\mbox{bc}(G)$ are isomorphic. The square of every such graph $G'$ does not contain a hamiltonian cycle (it is not possible to find a hamiltonian cycle in the square containing the third 2-valent vertex different from $i$, $j$, a contradiction), e.g. see the graph in Figure \ref{Counter} d), where $H_1,...,H_4$ are arbitrary connected graphs with at least one edge each. Note that conditions 4), 6) and the first part of condition 5) may hold.

\medskip

\noindent
\emph{Condition 6) does not hold.}
  
Hence $\sum_{t=1}^{n}m_{i}(B_t)<2k_i+\mbox{bn}(i)-2$ for some $i$ and consequently $m_i(B_t)=1$ for at least $3-\mbox{bn}(i)$ 2-blocks containing $i$. Note that, clearly, $\mbox{bn}(i)<2$ with respect to condition 3).

Let $r$ be the number of 2-blocks with $m_i(B_t)=1$. Each of these 2-blocks contains either exactly 2 cutvertices of $G$ or at least 3 cutvertices of $G$. Note that for 2-blocks containing only cutvertex $i$ we have $m_i(B_t)=2$ (see Algorithm 1 case d) below). We exchange every 2-block containing exactly 2 cutvertices of $G$ with a cycle $C_3$ and every 2-block containing $k$ cutvertices of $G$, $k\geq 3$, with a cycle $C_k$. In the first case note that we assume without loss of generality that there is no labelling such that we switch values 1 and 2 for both cutvertices of this 2-block to get a permissible labelling (again see Algorithm~1 case e) below).

Since $r\geq 3-\mbox{bn}(i)$, by the exchanging 2-block mentioned above we get a class of graphs $G'$ such that $\mbox{bc}(G')$ and $\mbox{bc}(G)$ are isomorphic. The square of every such graph $G'$ does not contain a hamiltonian cycle (if we try to construct a hamiltonian cycle in the square, then the degree of the cutvertex $i$ is at least 3, a contradiction), e.g. see graphs in Figure \ref{Counter} $e_1$) and $e_2$). For the graph in Figure \ref{Counter} $e_1$) it holds that $r=3-\mbox{bn}(i)=3-1=2$, the 2-block $B_1$ has exactly 2 cutvertices of $G$, the 2-block $B_2$ has $k=3$ cutvertices of $G$ (and hence $B_1$, $B_2$ are isomorphic to $C_3$) and $H_1,...,H_5$ are arbitrary connected graphs with at least one edge. For the graph in Figure \ref{Counter} $e_2$) it holds that $r=3-\mbox{bn}(i)=3-0=3$, the 2-block $B_1$ has exactly 2 cutvertices of $G$, the 2-block $B_2$ has $k=3$ cutvertices of $G$, the 2-block $B_3$ has $k=4$ cutvertices of $G$ (hence $B_1$, $B_2$ are isomorphic to $C_3$ and $B_3$ is isomorphic to $C_4$) and $H_1,...,H_7$ are arbitrary connected graphs with at least one edge. Note that conditions 4) and 5) may hold.

This finishes the proof of Theorem \ref{hamiltonian}.
\end{proof}

If there is a graph $G$ such that every labelling $m_i(B_t)$ violates at least one of the conditions 4) - 6) of Theorem \ref{hamiltonian}, then there is a graph $G'$ with $\mbox{bc}(G')=\mbox{bc}(G)$ such that $(G')^2$ is not hamiltonian  as it has been shown in the proof of Theorem~\ref{hamiltonian}. On the other hand, if we are able to construct a labelling $m_i(B_t)$ satisfying conditions 1) - 6) using the following algorithm, then $G^2$ is hamiltonian as it has been shown in the proof of Theorem~\ref{hamiltonian}.

\bigskip

\noindent
\emph{ALGORITHM 1:} 

Set $P_0=G-\cup_{t=1}^nB_t$. If any component of $P_0$ is not a caterpillar, then $\mbox{bn}(i)\geq 3$ for some $i\in\{1,2,...,s\}$ contradicting condition 4) in Theorem \ref{hamiltonian} and $G^2$ is not hamiltonian (e.g. see Figure \ref{Counter} a)). STOP.

If $G=P_0$, then $G$ is a caterpillar,  $n=0$ and $G^2$ is hamiltonian (see Observation in the proof of Theorem \ref{hamiltonian}) and $m_i(B_t)$ is not defined ($n=0$). STOP.

If $G$ is a 2-block, $G^2$ is hamiltonian by Theorem \ref{2-blockcycle} and $m_i(B_t)$ is not defined ($s=0$ and $n=1$). STOP.

We set $G_0=G-P_0$ and $m_i(B_t)=0$ if $i\notin V(B_t)$ for $i\in\{1,2,...,s\}$ and $t\in\{1,2,...,n\}$.

\bigskip

\noindent
START

We choose a 2-block $B$ containing at most 1 cutvertex of $G_0$. Note that $B$ is either a component of $G_0$ or an endblock of some component of $G_0$. If such endblock does not exist, we choose 2-block $B$ as a component of $G_0-H$ or an endblock of $G_0-H$ where $H$ is the union of all 2-blocks for which the labelling $m_i(B_t)$ is already set.  Let $c_1,c_2,...,c_k$ be all cutvertices of $G$ contained in $B$, $k\geq 1$. 

\medskip

\begin{itemize}
 \item[a)] If $k\geq 5$, then by condition 2) $m_{c_i}(B)\geq 1$ for $i=1,2,...,k$. Hence condition 5) in Theorem \ref{hamiltonian} does not hold and $G^2$ may not be hamiltonian (e.g. see Figure \ref{Counter} b)). STOP.

 \item[b)] If $k\geq 3$ and $\mbox{bn}(c_i)=2$ for some $i\in\{1,2,...,k\}$ , then by condition 3) $m_{c_i}(B)=2$ and by 2) $m_{c_j}(B)\geq 1$ for $j=1,2,...,k$. Hence condition~5) in Theorem \ref{hamiltonian} does not hold and $G^2$ may not be hamiltonian (e.g. see Figure \ref{Counter} c)). STOP.

 \item[c)] If $k=2$ and $\mbox{bn}(c_1)=\mbox{bn}(c_2)=2$, then by condition 3) $m_{c_1}(B)=2$ and $m_{c_2}(B)=2$. Hence condition 5) in Theorem \ref{hamiltonian} does not hold and $G^2$ may not be hamiltonian (e.g. see Figure \ref{Counter} d)). STOP.

 \item[d)] If $k=1$, then we set $m_{c_1}(B)=2$ (we maximize values $m_i(B_t)$ with respect to condition 6) in Theorem \ref{hamiltonian}). Note that, if the labelling $m_{i}(B_t)$ is set for all 2-blocks incident with $c_1$, then condition 6) holds for cutvertex $c_1$ with respect to the choice of $B$. 
 
 If the labelling $m_{i}(B_t)$ is set for all 2-blocks of $G$, then the labelling $m_i(B_t)$ satisfies the conditions of Theorem \ref{hamiltonian} and $G^2$ is hamiltonian. STOP.
 
 Otherwise we go to START. 

 \item[e)] If $k=2$ and $\mbox{bn}(c_i)\leq 1$ for $i\in\{1,2\}$, then we set $m_{c_1}(B)$ and $m_{c_2}(B)$ in the following way (without loss of generality $i=1$).
 
 Let $\mbox{bn}(c_2)=2$. Then we set $m_{c_1}(B)=1$ and $m_{c_2}(B)=2$ with respect to conditions 2), 3) and 5). 
 
 Let $\mbox{bn}(c_2)\leq 1$. Then for at least one of $c_1$, $c_2$ it holds that $m_{c_j}(B_t)$ for $j\in\{1,2\}$ is set for all 2-blocks $B_t$ except $B$ with respect to the choice of $B$ (again without loss of generality $j=1$). We set $m_{c_1}(B)=1$ and we verify condition 6) for $c_1$. If it holds, then we set $m_{c_2}(B)=2$ (again we maximize values $m_i(B_t)$ with respect to condition 6)). If condition~6) for $c_1$ does not hold for $m_{c_1}(B)=1$, then we set $m_{c_1}(B)=2$ and $m_{c_2}(B)=1$. 
 
 Now in both cases we verify condition 6) for $c_1$ and $c_2$ if the labelling $m_{c_1}(B_t)$ and $m_{c_2}(B_t)$ is set for all 2-blocks $B_t$. 
 
 If condition 6) does not hold in at least one case, then $G^2$ may not be hamiltonian (e.g. see Figure \ref{Counter} $e_1)$). STOP. 
 
 Hence suppose that condition 6) holds for $c_1$, $c_2$ if $m_{c_1}(B_t)$, $m_{c_2}(B_t)$ is set for all $B_t$, respectively.
 
 If the labelling $m_{i}(B_t)$ is set for all 2-blocks, then the labelling $m_i(B_t)$ satisfies the conditions of Theorem \ref{hamiltonian} and $G^2$ is hamiltonian. STOP.
  
 Otherwise we go to START. 
 
  \item[f)] If $k\in\{3,4\}$ and $\mbox{bn}(c_i)\leq 1$, then we set $m_{c_i}(B)=1$ for $i=1,2,...,k$. We verify condition 6) for all $c_i$ if the labelling $m_{c_i}(B_t)$ is set for all 2-blocks $B_t$. 
 
 If condition 6) does not hold in at least one case, then $G^2$ may not be hamiltonian (e.g. see Figure \ref{Counter} $e_2)$). STOP.
 
 Hence suppose that condition 6) holds for all $c_i$, $i=1,2,...,k$, for which $m_{c_i}(B_t)$ is set for all $B_t$.
 
 If the labelling $m_{i}(B_t)$ is set for all 2-blocks, then the labelling $m_i(B_t)$ satisfies the conditions of Theorem \ref{hamiltonian} and $G^2$ is hamiltonian. STOP.
 
 Otherwise we go to START. 
\end{itemize}

\noindent
PROOF OF THEOREM \ref{hamconnected}

\begin{proof}
 Let $x,y\in V(G)$. First we prove that there exists an $xy$-hamiltonian path $P$ in $G^2$ if there is no nontrivial bridge of $G$ and every block contains at most 2 cutvertices.

\bigskip
 
 (A) Suppose that $x$ and $y$ are in the same block $B$ of $G$. We proceed by induction on $n$, where $n$ is the number of blocks of $G$, $n\geq 1$.  
 
 For $n=1$, clearly $G=B$. If $B=K_2=xy$, then $G$ is also the $xy$-hamiltonian path in $G^2$ as required. If $B$ is a 2-block, then by Theorem \ref{strongF_3}, $G^2=B^2$ contains an $xy$-hamiltonian path $P$ as required.

\bigskip

Now suppose that the statement of Theorem \ref{hamconnected} is true for every graph with $n$ blocks and $G$ is a graph with $n+1$ blocks, $n\geq 1$. We distinguish 2 cases.

\begin{itemize}
 \item $B$ has exactly one cutvertex $c$.
 
 Without loss of generality we assume that $x\neq c$. If $B$ is a 2-block, then by Theorem \ref{strongF_3}, $B^2$ contains an $xy$-hamiltonian path $P_B$ containing an edge $cy'$ where $y'$ is a neighbor of $c$ in $B$. Note that $y'=x$ or $c=y$ is possible. If $B=K_2$, then $B=xy=y'c$ and $P_B=xy$ is an $xc$-hamiltonian path in $B^2$. By the induction hypothesis $(G-B)^2$ contains a $cc'$-hamiltonian path $P_G$ where $c'$ is a neighbor of $c$ in $G-B$. Then $P=P_B\cup P_G-cy'+y'c'$ is an $xy$-hamiltonian path in $G^2$ as required.
 
 \item $B$ has two cutvertices $c_1$, $c_2$.
 
 We denote by $G_1$, $G_2$ the two components of $G-B$ such that $c_i\in V(G_i)$ and let $c'_i$ be a neighbor of $c_i$ in $G_i$, $i=1,2$. By the induction hypothesis $(G_i)^2$ contains a $c_ic'_i$-hamiltonian path $P_{G_i}$, $i=1,2$.
 
 \begin{itemize}
  \item[a)] $c_i\notin \{x,y\}$ ($x$ and $y$ are not cutvertices).

 By Theorem \ref{F_4}, $B^2$ contains an $xy$-hamiltonian path $P_B$ containing the edges $c_iz_i$ where $z_i$ is a neighbor of $c_i$ in $B$, $i=1,2$. Note that $z_i\in\{x,y\}$ is possible. 
 
  \item[b)] Up to symmetry $c_1=x$ and $c_2\neq y$ (either $x$ or $y$ is a cutvertex of $G$).
  
  By Theorem \ref{strongF_3}, $B^2$ contains an $xy$-hamiltonian path $P_B$ containing the edges $c_iz_i$ where $z_i$ is a neighbor of $c_i$ in $B$, $i=1,2$. Note that $z_1=c_2$ or $z_2=y$ is possible. 
 
 \item[c)] $c_1=x$ and $c_2=y$ (similarly $c_1=y$ and $c_2=x$).
 
 By Theorem \ref{strongF_3ends}, $B^2$ contains an $xy$-hamiltonian path $P_B$ containing either the edges $c_iz_i$ where $z_i$ is a neighbor of $c_i$ in $B$, $i=1,2$, or the edges $c_1z_1$, $uv$ where $z_1$ is a neighbor of $c_1$ in $B$ and $u,v$ are neighbors of $c_2$ in $B$. 
 \end{itemize}
 
 In all cases except the case c), if $uv$ is the edge of $P_B$, $$P=P_{G_1}\cup P_B\cup P_{G_2}-\{c_1z_1,c_2z_2\}\cup\{c'_1z_1,c'_2z_2\}$$ is an $xy$-hamiltonian path in $G^2$ as required. 
 
 It remains to find an $xy$-hamiltonian path in $G^2$ if $uv$ is the edge of $P_B$. 
 
 If $G_2=K_2=c_2c'_2$, then $$P=P_{G_1}\cup P_B-\{c_1z_1,uv,c_2c'_2\}\cup\{c'_1z_1,c'_2u, c'_2v\}$$ is an $xy$-hamiltonian path in $G^2$ as required.    
 
 If $G_2\neq K_2$, then we prove that $(G_2)^2$ contains a hamiltonian cycle $C$ containing edges $c_2u_2$, $c_2v_2$ of $G_2$. Let $B_1,B_2,...,B_k$ be all 2-blocks of $G_2$ containing $c_2$. By Theorem \ref{2-blockcycle}, for  $i=1,2,...,k$, $(B_i)^2$ contains a hamiltonian cycle $C'_i$ containing three different edges $c_2u_2^i$, $c_2v_2^i$, $y_iy'_i$ of $B_i$ where $y_i$ is the second cutvertex of $G_2$ in $B_i$ if it exists.
 
 If $y_i$ exists, then we denote by $H_i$ a component of $G_2-(B_i-y_i)$ containing $y_i$. By the induction hypothesis $(H_i)^2$ contains a $y_id_i$-hamiltonian path $P_i$ where $d_i$ is a neighbor of $y_i$ in $H_i$. Then we set $C_i=C'_i\cup P_i-y_iy'_i+y'_id_i$. If $y_i$ does not exist, then we set $C_i=C'_i$.
 
 Let $T$ be the set of all leaves of $G_2$ adjacent to $c_2$. Then we find a cycle $C$ on $\cup_{i=1}^{k}V(C_i)\cup T$ by appropriately replacing edges $c_2u_2^i$, $c_2v_2^i$ with edges of $G^2$ joining $u_2^i$, $v_2^i$ in different $C_i$ and leaves adjacent to $c_2$ (similarly as in the proof of Theorem \ref{hamiltonian}) such that we preserve two edges ($c_2u_2^i$, $c_2v_2^i$ or $c_2l_1$, $c_2l_2$ where $l_1,l_2$ are two leaves of $G_2$ adjacent to $c_2$) as $c_2u_2$, $c_2v_2$. 
  
 Now $$P=P_{G_1}\cup P_B\cup C-\{c_1z_1,uv,c_2u_2,c_2v_2\}\cup\{c'_1z_1,u_2u, v_2v\}$$ is an $xy$-hamiltonian path in $G^2$ as required.

\end{itemize}

\bigskip 
 
 (B) Suppose that $x$ and $y$ are in different blocks of $G$.

 Let $P_G$ be any $xy$-path in $G$ and $c\in V(P_G)\setminus\{x,y\}$ be a cutvertex of $G$. Let $K$ be the component of $G-c$ containing $x$, $G_y=G-V(K)$ and $G_x=G-G_y$. Clearly $G_x\cup G_y=G$ and $G_x\cap G_y=c$.  If $G_x$, $G_y$ are isomorphic to $K_2$, then we set $P_x=G_x$, $P_y=G_y$, respectively. If $G_x$, $G_y$ are 2-blocks, then $(G_x)^2$, $(G_y)^2$ contains an $xc$-hamiltonian path $P_x$, a $cy$-hamiltonian path $P_y$ by Theorem \ref{strongF_3}, respectively. We proceed by induction on $n$, where $n$ is the number of blocks of $G$, $n\geq 2$.
  
 First assume that $G$ has exactly 2 blocks. Hence $G_x$, $G_y$ are isomorphic to $K_2$ or 2-blocks and $P=P_x\cup P_y$ is an $xy$-hamiltonian path in $G^2$ as required.
 
 Now suppose that the statement of Theorem \ref{hamconnected} is true for every graph with $n$ blocks and $G$ is a graph with $n+1$ blocks, $n\geq 2$. If $G_x$, $G_y$ is not a block, then by the induction hypothesis $(G_x)^2$, $(G_y)^2$ contains an $xc$-hamiltonian path $P_x$, a $cy$-hamiltonian path $P_y$, respectively. Then $P=P_x\cup P_y$ is an $xy$-hamiltonian path in $G^2$ as required.
  
\bigskip 

 Now it remains to prove that if there is a nontrivial bridge of $G$, then  $G^2$ is not hamiltonian connected and if $G$ contains a block containing more than 2 cutvertices, then there is a graph $G'$ such that $\mbox{bc}(G)$ and $\mbox{bc}(G')$ are isomorphic but $(G')^2$ is not hamiltonian connected. 
 
 Clearly, if there exists a nontrivial bridge $xy$ in $G$, then there is no $xy$-hamiltonian path in $G^2$ and $G^2$ is not hamiltonian connected.
 
 \begin{figure}[ht]
\begin{center}
\includegraphics[width=3cm]{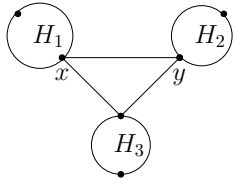}
\end{center}\caption{Graphs without $xy$-hamiltonian path in the square}\label{Counter2}
\end{figure} 

 Finally assume that $G$ contains a block $B$ containing $r$ cutvertices, where $r>2$. Then we exhange $B$ with a cycle $C_r$ and all other blocks with arbitrary blocks to get a class of graphs $G'$ such that $\mbox{bc}(G')$ and $\mbox{bc}(G)$ are isomorphic. Clearly the square of every such graph $G'$ does not contain a hamiltonian path between arbitrary two cutvertices of $G'$ in $C_r$ and hence $(G')^2$ is not hamiltonian connected, e.g. with Figure \ref{Counter2}, where $r=3$ and $H_1,H_2,H_3$ are arbitrary connected graphs with at least one edge.
\end{proof}

Similarly as for Theorem \ref{hamiltonian} we state the following algorithm to verify conditions of Theorem \ref{hamconnected}.

\bigskip

\noindent
\emph{ALGORITHM 2:} 

Let $G'=G-S$ where $S$ is the set of all endblocks of $G$. Let $\mbox{cvn}_G(B)$ be the number of cutvertices of $G$ in $B$.

\noindent
START

Find an endblock $B$ of $G'$.

\begin{itemize}
 \item If $B$ is a bridge of $G'$, then $B$ is a nontrivial bridge of $G$ and $G^2$ is not hamiltonian connected. STOP.
 \item Let $B$ be a 2-block.  
   \begin{itemize}
     \item If $\mbox{cvn}_G(B)>2$, then $G^2$ may not be hamiltonian connected (e.g. see Figure \ref{Counter2}). STOP.
     \item If $\mbox{cvn}_G(B)\leq 2$, then $G':=G'-B$.
        \begin{itemize}
         \item If $G'=\emptyset$, then $G^2$ is hamiltonian connected. STOP.
         \item If $G'\neq\emptyset$, then go to START.
        \end{itemize}
   \end{itemize}
\end{itemize}

In both algorithms in this paper, determining blocks and especially endblocks and bridges, cutvertices, block-cutvertex graphs, and the parameters $\mbox{bn}(i)$, $\mbox{cvn}_G(B)$ can be determined in polynomial time. 

As a consequence, polynomial running time in Algorithm 2 is guaranteed. For, determining (potentially) not being Hamiltonian connected, can be determined instantly once a nontrivial bridge, a block with more than 2 cutvertices has been found. And deleting an endblock reduces the size of $G'$ linearly. 

Now consider the running time of Algorithm 1. The first decision to be made is whether $P_0$  is a forest of caterpillars – this can be done in linear time. After that, at every step 'one chooses a 2-block $B$ as a component of $G_0-H$ or an endblock of $G_0-H$ where $H$ is the union of all 2-blocks for which the labelling $m_i(B_t)$ is already set'. Clearly, identifying such $B$ can be done in linear time. The same applies to working through the cases for defining the various values of $m_i(B)$. 

Summarizing, it follows that both algorithms run in polynomial time. We note however, that these algorithms can only decide the existence or potential non-existence of hamiltonian cycles or hamiltonian paths in the square of graphs under consideration; they do not construct any such cycle or path.

\section{Conclusion}
The main results of this paper are Theorem \ref{hamiltonian} and Theorem \ref{hamconnected}.
As we mention in Introduction Fleischner in \cite{Fle} proved that in the square of graphs hamiltonicity and vertex-pancyclicity are equivalent concepts, and so are hamiltonian connectedness and panconnectedness. Hence we proved in fact that for graphs satisfying assumptions of Theorem \ref{hamiltonian}, Theorem \ref{hamconnected} the square of these graphs is vertex-pancyclic, panconnected, respectively.
\medskip

As an easy corollary of Theorem \ref{hamconnected}  we get the following result.

\begin{corollary}
\label{Block-chain}
 Let $G$ be a block-chain. Then $G^2$ is panconnected if and only if every innerblock of $G$ is a 2-block. 
\end{corollary}

Moreover Corollary \ref{Block-chain} is also the answer to Problem 1 stated by Chia et al. in \cite{ChiaOngTan} that for a graph $G$ with only two cutvertices it is true that $G^2$ is panconnected if and only if the unique block containing the two cutvertices is not the complete graph on two vertices.

\medskip

\noindent {\bf Acknowledgements}.
This work was partly supported by the European Regional Development Fund (ERDF), project NTIS - New Technologies for Information Society, European Centre of Excellence, CZ.1.05/1.1.00/02.0090.

The first author was partly supported by project GA20-09525S of the Czech Science Foundation. The second author was supported in part by FWF grant P27615.

\end{document}